# $M|G|\infty$ QUEUE BUSY PERIOD WITH PME DISTRIBUTION[1]

Prof. Dr. **MANUEL ALBERTO M. FERREIRA**
Instituto Universitário de Lisboa (ISCTE – IUL), BRU - IUL, Lisboa, Portugal
manuel.ferreira@iscte.pt

**ABSTRACT**

In this article it is shown that if the busy period of a $M|G|\infty$ queue system is PME distributed, the service time is a random variable with a long-tail distribution. The result is obtained through Laplace transforms analysis.

## 1. INTRODUCTION

In a $M|G|\infty$ queue system, $\lambda$ is the Poisson process arrivals rate, $\alpha$ is the mean service time, $G(.)$ is the service time distribution function and so $\alpha = \int_0^\infty [1 - G(t)]\, dt$. The traffic intensity is $\rho = \lambda\alpha$ and $B$ is the busy period length.

Note the importance of the busy period study, for this queuing system, because in it any customer, when it arrives, finds immediately an available server. So the problem is "for how long the servers – and how many servers? – must be available? That is: how long is the busy period length?"

It will be supposed that the $B$ probability density function is

$$g_r(t) = \int_{\frac{r-1}{r}}^{\infty} f_r(y) y^{-1} e^{-\frac{t}{y}} dy, r > 1 \qquad (1.1)$$

where $f_r(y) = r\left(\frac{r-1}{r}\right)^r x^{-(r+1)}, x \geq \frac{r-1}{r}$, is a Pareto distribution probability density function. It is a PME[2], see (1), and it is long-tail type. The $g_r$ moments are

$$m_n = n! \frac{r}{r-n}\left[\frac{r-1}{r}\right]^n, n = 1,2,\ldots \qquad (1.2).$$

---

[1] This work was financially supported by FCT through the Strategic Project PEst-OE/EGE/UI0315/2011.

[2] **P**areto **M**ixture of **E**xponentials.

Through Laplace transform analysis it will be emphasized that in these circumstances the service time is a random variable with a long-tail distribution.

## 2. THE PME LAPLACE TRANSFORM

Calling $\hat{g}_r(s)$ the Laplace transform of a PME with parameter $r$, see again (1),

$$\hat{g}_r(s) = r\left(\frac{r-1}{r}\right)^r \int_0^{\frac{r}{r-1}} \frac{x^r}{s+x} dx, r > 1 \qquad (2.1).$$

But $\hat{g}_r(s) = r\left(\frac{r-1}{r}\right)^r \left(\int_0^{\frac{r}{r-1}}\left(x^{r-1} - \frac{sx^{r-1}}{s+x}\right)dx\right) = r\left(\frac{r-1}{r}\right)^r \left(\int_0^{\frac{r}{r-1}} x^{r-1} dx - s\int_0^{\frac{r}{r-1}} \frac{x^{r-1}}{s+x} dx\right) = r\left(\frac{r-1}{r}\right)^r \left(\left[\frac{x^r}{r}\right]_0^{\frac{r}{r-1}} - s\int_0^{\frac{r}{r-1}} \frac{x^{r-1}}{s+x} dx\right) = r\left(\frac{r-1}{r}\right)^r \left(\frac{1}{r}\left(\frac{r}{r-1}\right)^r - s\int_0^{\frac{r}{r-1}} \frac{x^{r-1}}{s+x} dx\right) = 1 - sr\left(\frac{r-1}{r}\right)^r \int_0^{\frac{r}{r-1}} \frac{x^{r-1}}{s+x} dx.$

If $\hat{h}_r(s)$ is the PME with parameter $r$ tail Laplace transform, as $\hat{h}_r(s) = \frac{1}{s} - \frac{1}{s}\hat{g}_r(s)$,

$$\hat{h}_r(s) = r\left(\frac{r-1}{r}\right)^r \int_0^{\frac{r}{r-1}} \frac{x^{r-1}}{s+x} dx, r > 1 \quad (2.2).$$

So, after (2.2),

$$\hat{h}_r^{(n)}(s) = (-1)^n n! \, r\left(\frac{r-1}{r}\right)^r \int_0^{\frac{r}{r-1}} \frac{x^{r-1}}{(s+x)^{n+1}} dx, n = 0,1,2,\ldots \qquad (2.3)$$

where $\hat{h}_r^{(n)}$ is the n$^{th}$ order derivative of $\hat{h}_r$.

Then $\hat{h}_r^{(n)}(0) = (-1)^n n! \, r\left(\frac{r-1}{r}\right)^r \int_0^{\frac{r}{r-1}} x^{r-n-2} dx, n = 0,1,2,\ldots$ . But $\int_0^{\frac{r}{r-1}} x^{r-n-2} dx = \begin{cases} \left[\frac{x^{r-n-1}}{r-n-1}\right]_0^{\frac{r}{r-1}}, n \neq r-1 \\ [\log|x|]_0^{\frac{r}{r-1}}, n = r-1 \end{cases} = \begin{cases} \frac{\left(\frac{r}{r-1}\right)^{r-(n+1)}}{r-(n+1)}, r > n+1 \\ -\infty, r \leq n+1 \end{cases}$.

So $\hat{h}_r^{(n)}(0) = \begin{cases} (-1)^n n! \, r\left(\frac{r-1}{r}\right)^r \frac{\left(\frac{r}{r-1}\right)^{r-(n+1)}}{r-(n+1)}, r > n+1 \\ (-1)^n(-\infty), 1 < r \leq n+1 \end{cases}$ or, equivalently,

$$\hat{h}_r^{(n)}(0) = \begin{cases} (-1)^n n! \, r \left(\dfrac{r-1}{r}\right)^r \dfrac{\left(\dfrac{r}{r-1}\right)^{r-(n+1)}}{r-(n+1)}, & n < r+1, r > 1 \\ (-1)^n(-\infty), & n \geq r-1 \end{cases} \quad (2.4).$$

## 3. $M|G|\infty$ BUSY PERIOD TAIL LAPLACE TRANSFORM

Call $U(t)$ the $M|G|\infty$ busy period tail and $u(s)$ the respective Laplace transform so, see (2),

$$\frac{d}{dt}\left(\frac{1 - e^{-\lambda \int_0^t [1-G(v)]dv}}{1 - e^{-\rho}}\right) = TL^{-1}\left(\frac{1}{1-e^{-\rho}} \frac{\lambda u(s)}{\lambda u(s)+1}\right) \quad (3.1)$$

and

$$\frac{e^{-\lambda \int_0^t [1-G(v)]dv} \lambda(1-G(t))}{1-e^{-\rho}} = TL^{-1}\left(\frac{1}{1-e^{-\rho}} \frac{\lambda u(s)}{\lambda u(s)+1}\right) \quad (3.2)$$

being $TL^{-1}(.)$ the inverse Laplace transform.

So, $\int_0^\infty t^n \dfrac{e^{-\lambda \int_0^t [1-G(v)]dv} \lambda(1-G(t))}{1-e^{-\rho}} dt =$
$(-1)^n \dfrac{1}{1-e^{-\rho}} \left(\dfrac{\lambda u(s)}{\lambda u(s)+1}\right)^{(n)}_{s=0} \Leftrightarrow \int_0^\infty t^n e^{-\lambda \int_0^t [1-G(v)]dv} \lambda(1-G(t)) dt =$
$(-1)^n \left(\dfrac{\lambda u(s)}{\lambda u(s)+1}\right)^{(n)}_{s=0}.$

As $e^{-\lambda \int_0^t [1-G(v)]dv} \leq 1$ the consequence is that $\int_0^\infty t^n \lambda (1-G(t)) dt \geq$
$(-1)^n \left(\dfrac{\lambda u(s)}{\lambda u(s)+1}\right)^{(n)}_{s=0} \Leftrightarrow \int_0^\infty \dfrac{1-G(t)}{\alpha} dt \geq \dfrac{(-1)^n}{\rho} \left(\dfrac{\lambda u(s)}{\lambda u(s)+1}\right)^{(n)}_{s=0}$ this entire happening if
$\dfrac{1}{1-e^{-\rho}} \dfrac{\lambda u(s)}{\lambda u(s)+1}$ is in fact a probability density function Laplace transform, being enough so that this happens that $(-1)^n \left(\dfrac{\lambda u(s)}{\lambda u(s)+1}\right)^{(n)}_{s=0} \geq 0, n = 0,1,2,\ldots$.

## 4. $M|G|\infty$ BUSY PERIOD WITH PME DISTRIBUTION

Note that in the $\frac{\lambda u(s)}{\lambda u(s)+1}$ n$^{th}$ order derivative, $u^{(n)}(s)$ always appears with a positive sign in the numerator and, for $s = 0$, if $u(s)$ is given by (2.1), in that order n derivative the denominator is $(\lambda + 1)^{n+1}$. It is enough to take in account the quotient derivative expression and that $\hat{h}_r(0) = 1$.

So, after (2.4), it is concluded that

-If there is a service distribution such that the $M|G|\infty$ queue busy period is distributed as a PME distribution with parameter $r$, the service equilibrium distribution moments of order greater than $r - 1$, centered in the origin, are infinite.

Note that

-The service equilibrium distribution, with these moments, is a long-tail distribution, see again (1),

-As $\frac{u(s)}{\frac{e^\rho-1}{\lambda}}$ is the $M|G|\infty$ queue busy period equilibrium distribution Laplace transform it is also concluded that if this has moments of order greater than $n, n \in \mathbb{N}$, infinite, the same happens with the service equilibrium distribution, that is: they are both long –tail distributions.

## 5. CONCLUDING REMARKS

As it is stated in (1) the PME are long-tail distributions. So it is checked in this work that there is an uncontested association between long –tail service distributions and long-tail busy period distributions for the $M|G|\infty$ queue, as it was shown in (3).

The PME were introduced in (1). There they were a tool to investigate properties of waiting times tail probabilities in queues with long-tail service-time distributions. For this investigation the authors developed algorithms for computing the waiting time distribution by Laplace transform inversion when the Laplace transforms of the inter-arrival time and service time distributions are known. The procedure here followed is similar.